\input ./style/arxiv-general.cfg
\documentclass[MSNbibl,number,citesort,seceqn,dvips]{arxbj}
\makeatletter
   \@ifpackageloaded{graphicx}{}{\usepackage{graphicx}}
\makeatother
\usepackage{mathbh}

\volume{23}
\issue{1}
\pubyear{2017}
\firstpage{405}
\lastpage{431}
\doi{10.3150/15-BEJ747}
\docsubty{FLA}

\makeatletter
\newcommand{\rrvert}{\vert}
\newcommand{\llvert}{\vert}
\renewcommand{\mid}{|}
\newtheorem{theorem}{Theorem}[section]
\newtheorem{conjecture}[theorem]{Conjecture}
\newproclaim{definition}[theorem]{Definition}
\newtheorem{lemma}[theorem]{Lemma}
\newtheorem{proposition}[theorem]{Proposition}
\newremark{remark}[theorem]{Remark}

\newcommand{\one}{\mathbh{1}}
\newcommand{\nn}{\nonumber}
\newcommand{\n}{\mathbf{n}}

\def\bE{\mathbb{E}}
\def\bN{\mathbb{N}}
\def\bP{\mathbb{P}}
\def\bR{\mathbb{R}}
\def\bZ{\mathbb{Z}}
\def\cG{\mathcal{G}}
\def\cK{\mathcal{K}}
\def\cL{\mathcal{L}}
\def\cM{\mathcal{M}}
\def\cR{\mathcal{R}}
\def\cZ{\mathcal{Z}}
\def\kI{\mathfrak{I}}
\def\kS{\mathfrak{S}}
\def\s{\sigma}
\def\w{\omega}
\def\e{\varepsilon}
\def\essinf{\mathop{\operatorname{ess}\operatorname{inf}}}
\def\esssup{\mathop{\operatorname{ess}\operatorname{sup}}}
\makeatother

\begin{document}
\begin{frontmatter}

\title{Variational formulas and disorder regimes of random walks in random potentials}
\runtitle{Variational formulas and disorder regimes of RWRP}

\begin{aug}
\author[A]{\inits{F.}\fnms{Firas}~\snm{Rassoul-Agha}\thanksref{A}\ead[label=e1]{firas@math.utah.edu}\ead[label=u1,url]{http://www.math.utah.edu/\textasciitilde firas}},
\author[B]{\inits{T.}\fnms{Timo}~\snm{Sepp\"al\"ainen}\thanksref{B}\ead[label=e2]{seppalai@math.wisc.edu}\ead[label=u2,url]{http://www.math.wisc.edu/\textasciitilde seppalai}}
\and
\author[C]{\inits{A.}\fnms{Atilla}~\snm{Yilmaz}\corref{}\thanksref{C}\ead[label=e3]{atillayilmaz@ku.edu.tr}\ead[label=u3,url]{http://home.ku.edu.tr/\textasciitilde atillayilmaz/}}
\address[A]{Department of Mathematics, University of Utah, 155S 1400E,
Salt Lake City, UT 84112, USA.\\ \printead{e1,u1}}
\address[B]{Department of Mathematics, University of
Wisconsin--Madison, Van Vleck Hall, 480 Lincoln Dr., Madison WI
53706-1388, USA.\\ \printead{e2,u2}}
\address[C]{Department of Mathematics, Ko\c{c} University, 34450 Sar{\i}yer, Istanbul, Turkey.\\ \printead{e3,u3}}
\end{aug}

%
\received{\smonth{10} \syear{2014}}
%
\revised{\smonth{4} \syear{2015}}

%
\begin{abstract}
We give two variational formulas \textup{(\ref{qVar1})} and \textup{(\ref{qVar2})} for the quenched
free energy of a random walk
in random potential (RWRP) when (i) the underlying walk is directed or
undirected, (ii) the environment
is stationary and ergodic, and (iii) the potential is allowed to depend
on the next step of the walk
which covers random walk in random environment (RWRE). In the directed
i.i.d. case, we also give two
variational formulas \textup{(\ref{aVar1})} and \textup{(\ref{aVar2})} for the annealed free energy
of RWRP. These four formulas are
the same except that they involve infima over different sets, and the
first two are modified versions
of a previously known variational formula \textup{(\ref{qVar0})}  for which we provide
a short alternative proof.
Then, we show that \textup{(\ref{qVar0})} always has a minimizer, \textup{(\ref{aVar2})} never has
any minimizers unless the
RWRP is an RWRE, and \textup{(\ref{aVar1})} has a minimizer if and only if the RWRP is
in the weak disorder
regime. In the latter case, the minimizer of \textup{(\ref{aVar1})} is unique and it
is also the unique minimizer of
\textup{(\ref{qVar1})}, but \textup{(\ref{qVar2})} has no minimizers except for RWRE. In the case of
strong disorder, we give a
sufficient condition for the nonexistence of minimizers of \textup{(\ref{qVar1})} and
\textup{(\ref{qVar2})} which is satisfied
for the log-gamma directed polymer with a sufficiently small parameter.
We end with a conjecture
which implies that \textup{(\ref{qVar1})} and \textup{(\ref{qVar2})} have no minimizers under very
strong disorder.
\end{abstract}

%
\end{frontmatter}


\section{Introduction}

\subsection{The model}\label{themodelsec}

Random walk in random potential (RWRP) on $\bZ^d$, with $d\ge1$, has
three ingredients.
\begin{longlist}[(ii)]

\item[(i) \textit{The underlying walk}:] Fix a finite set $\cR\subset
\bZ^d$ with $\llvert  \cR \rrvert  \ge2$. Define $p:\bZ^d\to[0,1]$ by $p(z) =
1/\llvert  \cR \rrvert  $ if $z\in\cR$ and $p(z) = 0$ otherwise. Consider random walk
on $\bZ^d$ with i.i.d. steps that have $p$ as their common
distribution. This walk induces a probability measure $P_x$ on paths
starting at $x\in\bZ^d$. Expectations under $P_x$ are denoted by $E_x$.

\item[(ii) \textit{The environment}:] Let $\cG$ be the additive subgroup
of $\bZ^d$ generated by $\cR$. Take a probability space $(\Omega,\kS
,\bP)$ equipped with an Abelian group $\{T_x: x\in\cG\}$ of
measurable transformations such that (i) $T_{x+y} = T_x\circ T_y$ and
(ii) $T_0$ is the identity. Assume that $\bP$ is invariant and ergodic\vadjust{\goodbreak}
w.r.t. this group. Expectations under $\bP$ are denoted by $\bE$,
and sample points from $(\Omega,\kS,\bP)$ are referred to as environments.

\item[(iii)] \textit{The potential}: Take a measurable function $V:\Omega
\times\cR\to\bR$. For every $\w\in\Omega$, $x\in\bZ^d$ and
$z\in\cR$, the quantity $V(T_x\w,z)$ is referred to as the potential
at the ordered pair $(x,x+z)$ in the environment $\w$. 
\end{longlist}
Given $n\ge1$ and $\w\in\Omega$, we define the quenched RWRP
probability measure
\[
Q_{n,x}^\w\bigl((X_i)_{i\ge0}\in\cdot
\bigr) = \frac{1}{\cZ_{n,x}^\w}E_x \bigl[e^{\sum_{i=0}^{n-1}V(T_{X_i}\w,Z_{i+1})}
\one_{\{(X_i)_{i\ge0}\in
\cdot\}} \bigr]
\]
on paths starting at any $x\in\bZ^d$. Here, $(X_i)_{i\ge0}$ denotes
the random path with increments $Z_{i+1} = X_{i+1} - X_i$, and
\[
\cZ_{n,x}^\w= E_x \bigl[e^{\sum_{i=0}^{n-1}V(T_{X_i}\w
,Z_{i+1})} \bigr]
\]
is the normalizing factor, called the quenched partition function.


\begin{remark}\label{gerektiyaw}
We have fixed $p$ to be the uniform distribution on $\cR$, but we can
easily incorporate more general cases. Indeed, consider a measurable
$\hat p:\Omega\times\bZ^d\to[0,1]$ such that, for $\bP$-a.e. $\w
$:
(i) $\hat p(\w,z) > 0$ if and only if $z\in\cR$; and (ii) $\sum_{z\in\cR}\hat p(\w,z) = 1$. Then, the discrete-time Markov chain on
$\bZ^d$, with transition probabilities $\pi_{x,y}^\w:= \hat p(T_x\w
,y-x)$ for $x,y\in\bZ^d$, is a quenched random walk in random
environment (RWRE). Taking the underlying walk to be this RWRE is
equivalent to adding $- \log\hat p(\w,z) - \log\llvert  \cR \rrvert  $ to the
potential $V(\w,z)$.
\end{remark}

\begin{remark}
We have given a rather abstract definition of the environment space.
The canonical setting is as follows: there is a Borel set $\Gamma
\subset\bR$, and $(\Omega,\kS)$ is $\Gamma^{\bZ^d}$ equipped with
the product Borel $\s$-algebra. In this case, environments are
represented as $\w= (\w_x)_{x\in\bZ^d}$, and elements of the group
$\{T_x: x\in\cG\}$ are translations defined by $(T_x\w)_y = \w_{x+y}$.
\end{remark}

In the initial parts of this paper, we will consider RWRP with the
abstract environment formulation. However, in later parts, we will
adopt the canonical model and make the following extra assumptions.
\begin{longlist}[(Dir)]

\item[(Dir) {\it Directed nearest-neighbor walk}:] $\cR= \{e_1,\ldots, e_d\}$, the standard basis for $\bR^d$, with \mbox{$d\ge2$}.

\item[(Ind) {\it Independent environment}:] The components of $\w=
(\w_x)_{x\in\bZ^d}$ are i.i.d. under $\bP$.

\item[(Loc) {\it Local potential}:] There exists a $V_o:\Gamma\times
\cR\to\bR$ such that $V(\w,z) = V_o(\w_0,z)$ for every $\w= (\w
_x)_{x\in\bZ^d} \in\Omega= \Gamma^{\bZ^d}$ and $z\in\cR$.
\end{longlist}
These assumptions enable us to use martingale techniques in the
analysis of the asymptotic behaviour of RWRP, see Section~\ref{directedsec}. If $V_o$ does not depend on $z$, then RWRP is also
referred to as a directed polymer. However, we prefer to keep the $z$
dependence because, this way, the results on the quenched free energy
of RWRP have implications regarding large deviations, see Remark~\ref{zelazim}.

There is a vast literature on RWRP, RWRE and directed polymers: see the
lectures/surveys \cite
{BolSzn2002,ComShiYos2004,Hol2009,Gia2007,Szn2004,Zei2004} and the
references therein. In what follows, we will focus only on the parts of
the literature that are directly relevant to our results.

%

\subsection{Quenched free energy and large deviations}\label{quenchedLDPsec}

In a recent paper \cite{RasSepYil2013}, we prove the $\bP$-a.s. existence of the quenched free energy
%
\begin{equation}
\label{freirg} \Lambda_q(V):= \lim_{n\to\infty}
\frac{1}{n}\log\cZ_{n,0}^\w.
\end{equation}
In order to give the precise statement of this result, we need two definitions.

\begin{definition}\label{classK}
A measurable function $F:\Omega\times\cR\to\bR$ is said to be a
centered cocycle if it satisfies the following conditions.
\begin{longlist}[(ii)]
\item[(i)] Centered: $\bE[\llvert  F(\cdot,z) \rrvert  ] < \infty$ and $\bE[F(\cdot
,z)] = 0$ for every $z\in\cR$.


\item[(ii)] Cocycle:
\[
\sum_{i=0}^{m-1}F(T_{x_i}
\w,z_{i+1}) = \sum_{j=0}^{n-1}F
\bigl(T_{x'_j}\w,z'_{j+1}\bigr)
\]
for $\bP$-a.e. $\w$, every $m,n\ge1$, $(x_i)_{i=0}^m$ and
$(x'_j)_{j=0}^n$ such that $z_{i+1}:= x_{i+1} - x_i\in\cR$,
$z'_{j+1}:= x'_{j+1} - x'_j\in\cR$, $x_0 = x'_0$ and $x_m = x'_n$.
\end{longlist}
The class of centered cocycles is denoted by $\cK_0$.
\end{definition}

\begin{definition}\label{classL}
A measurable function $V:\Omega\times\cR\to\bR$ is said to be in
class $\cL$ if $\bE[\llvert  V(\cdot,z) \rrvert  ] < \infty$ and 
\[
\limsup_{\delta\to0}\limsup_{n\to\infty}\max
_{x\in\bigcup
_{j=1}^nD_j}\frac{1}{n}\sum_{0\le i\le\delta n}
\bigl\llvert V(T_{x+iz'}\w,z)\bigr\rrvert = 0
\]
for $\bP$-a.e. $\w$ and every $z,z'\in\cR$ such that $z'\ne 0$, where
%
\begin{equation}
\label{mrdj} D_j = \bigl\{z_1 + \cdots+
z_j\in\bZ^d: z_i\in\cR\mbox{ for every }
i=1,\ldots,j\bigr\}
\end{equation}
denotes the set of points accessible from the origin in exactly $j$
steps chosen from $\cR$.
\end{definition}

\begin{theorem}\label{KVarOld}
Assume that $\kS$ is countably generated and $V\in\cL$. Then, the
limit in (\ref{freirg}) exists $\bP$-a.s., is deterministic, and satisfies
%
\begin{equation}
\label{bboduk} \Lambda_q(V) = \inf_{F\in\cK_0}\bP
\mbox{-}\esssup_{\w} \biggl\{ \log \biggl(\sum
_{z\in\cR}p(z)e^{V(\w,z) + F(\w,z)} \biggr) \biggr\}\in(-\infty,
\infty].
\end{equation}
\end{theorem}

This result was initially obtained in \cite{Yil2009b} for bounded
potentials under the assumption that $\{\pm e_1,\ldots,\pm e_d\}
\subset\cR$. The version in Theorem~\ref{KVarOld} is part of \cite{RasSepYil2013}, Theorem 2.3, which is valid for potentials of the form
$V:\Omega\times\cR^\ell\to\bR$ with arbitrary $\ell\ge1$.
Actually, the latter result contains two variational formulas for
$\Lambda_q(V)$, but the second one is not directly relevant for our
purposes in this paper, so we omit it for the sake of brevity.

The proof of Theorem~\ref{KVarOld} is based on a rather technical
approach involving careful applications of ergodic and minimax
theorems, which was developed in \cite{KosRezVar2006,KosVar2008} in
the context of stochastic homogenization of viscous Hamilton--Jacobi
equations and was first adapted in \cite{Ros2006} to large deviations
for RWRE. However, the existence of the a.s. limit in (\ref{freirg})
can be shown more easily (without giving any formulas for $\Lambda
_q(V)$) by subadditivity arguments and additional estimates such as
concentration inequalities or lattice animal bounds. This has been done
in \cite{CarHu2002,ComShiYos2003,Var2007} for directed polymers under
various moment assumptions on the potential, and more recently in \cite{RasSep2014}, Theorem 2.2(b), in the setting of Theorem~\ref{KVarOld}.
In fact, the latter result drops the assumption that $\kS$ is
countably generated and only requires $V\in\cL$. We record it below
for future reference.

\begin{theorem}\label{pe2pe}
Assume that $V\in\cL$. Then, the limit in (\ref{freirg}) exists $\bP
$-a.s., is deterministic, and satisfies $\Lambda_q(V)\in(-\infty
,\infty]$.
\end{theorem}

The hypotheses of Theorem~\ref{KVarOld} are satisfied in the commonly
studied examples. First of all, in the canonical setting, the product
Borel $\sigma$-algebra is countably generated. Second, bounded
potentials are in class $\cL$ under arbitrary stationary and ergodic
$\bP$, and so is any $V$ with $\bE[\llvert  V(\cdot,z) \rrvert  ]<\infty$ when
$d=1$. In the multidimensional case under \textup{(Ind)} and \textup{(Loc)}, it suffices
to have $\bE[\llvert  V(\cdot,z) \rrvert  ^p]<\infty$ for some $p>d$. In general,
there is a tradeoff between the degree of mixing in $\bP$ and the
moment of $V(\cdot,z)$ required. See \cite{RasSepYil2013}, Lemma A.4,
for further details and proofs. Note that these assumptions do not rule
out $\Lambda_q(V) = \infty$. Indeed, it is easy to see that the
latter holds under \textup{(Ind)} and \textup{(Loc)} when $\cR$ allows multiple visits
to points and $V$ is unbounded.

Assume additionally that $\Omega$ is a compact metric space and $\kS$
is its Borel $\s$-algebra. (These assumptions are valid in the
canonical setting if $\Gamma$ is compact.) Let $\cM_s(\Omega\times
\cR)$ be the space of Borel probability measures $\mu$ on $\Omega
\times\cR$ such that
\[
\sum_{z\in\cR}\int_\Omega\varphi(\w)
\mu(d\w,z) = \sum_{z\in
\cR}\int_\Omega
\varphi(T_z\w)\mu(d\w,z)
\]
for every $\varphi\in C_b(\Omega)$, where $C_b(\cdot)$ denotes the
space of bounded continuous functions. It is shown in \cite{RasSepYil2013}, Theorem~3.1, that, when $\Lambda_q(V) < \infty$, Theorem~\ref
{KVarOld} implies a large deviation principle (LDP) for the quenched
distributions $Q_{n,0}^\w(R_n\in\cdot)$ on $\cM_s(\Omega\times\cR
)$ of the empirical measure
\[
R_n = \frac{1}{n}\sum_{i=0}^{n-1}
\delta_{T_{X_i}\w,Z_{i+1}}.
\]
The rate function of this LDP has the following formula:
%
\begin{equation}
\label{level2RF} \kI_q(\mu) = \sup_{f\in C_b(\Omega\times\cR)} \biggl\{
\sum_{z\in
\cR}\int_{\Omega} f(\w,z)\mu(d
\w,z) - \Lambda_q(f+V) \biggr\} + \Lambda_q(V).
\end{equation}
As a corollary, we get an LDP for $Q_{n,0}^\w(X_n/n\in\cdot)$ with
the rate function
%
\begin{equation}
\label{level1RF} I_q(v) = \sup_{\lambda\in\bR^d} \bigl\{\lambda
\cdot v - \Lambda _q(f_\lambda+ V) \bigr\} +
\Lambda_q(V).
\end{equation}
Here, $f_\lambda:\Omega\times\cR\to\bR$ is defined by $f_\lambda
(\w,z) = \lambda\cdot z$.
%
\begin{remark}\label{zelazim}
Observe that $f_\lambda+ V$ depends on $z$ even if $V$ does not. This
is why it is important to allow potentials that depend on $z$ in
Theorem~\ref{KVarOld}.
\end{remark}
%
\begin{remark}\label{haybeyok}
Theorem~\ref{pe2pe} also implies the aforementioned LDPs, but does not
provide formulas for $\Lambda_q(V)$, $\Lambda_q(f+V)$ and $\Lambda
_q(f_\lambda+V)$ appearing in the rate functions (\ref{level2RF}) and
(\ref{level1RF}).
\end{remark}

In the theory of large deviations, the former LDP is referred to as
level-2, and gives the latter one (known as level-1) via the so-called
contraction principle. See \cite
{DemZei2010,Hol2000,DeuStr1989,RasSep2010,Var1984} for the definitions
of these concepts as well as general background on large deviations.
The highest level is level-3 (also known as process level) and is
established for RWRP in \cite{RasSepYil2013}, Theorem 3.2. This last
LDP covers and strengthens various previous results on the quenched
large deviations for RWRP and RWRE such as \cite
{GreHol1994,Zer1998a,Zer1998b,ComGanzei2000,Var2003,CarHu2004,Ros2006,Yil2009a,Yil2009b,AveHolRed2010,RasSep2011}.
See \cite{RasSepYil2013}, Section~1.3, for a detailed account.

%

\subsection{Directed i.i.d. case: Disorder regimes}\label{directedsec}

Assume that the conditions \textup{(Dir)}, \textup{(Ind)} and \textup{(Loc)} from Section~\ref{themodelsec} are satisfied. In this case, the $\s$-algebras
\begin{eqnarray*}
\kS_0^n &=& \s\bigl(\w_x: x\in
\bZ_+^d, \llvert x \rrvert _1 \le n-1\bigr)\quad\mbox{and}
\\
\kS_0^\infty&=& \s\bigl(\w_x: x\in
\bZ_+^d\bigr)
\end{eqnarray*}
are relevant. Here and throughout, $\bZ_+:= \bN\cup\{0\}$ is the
set of nonnegative integers and $\llvert  x \rrvert  _1:= \llvert  x_1\rrvert   +\cdots+\llvert  x_d\rrvert  $ denotes
the $\ell_1$-norm. Functions that are measurable w.r.t. $\kS
_0^\infty$ are sometimes referred to as future measurable.

Define the annealed free energy
%
\begin{equation}
\label{avratdef} \Lambda_a(V):= \log \biggl(\sum
_{z\in\cR}p(z)\bE \bigl[e^{V(\cdot
,z)} \bigr] \biggr)\in(-\infty,
\infty].
\end{equation}
It is straightforward to check that
\[
W_n(\w):= \frac{\cZ_{n,0}^\w}{\bE[\cZ_{n,0}^\w]} = E_0
\bigl[e^{\sum_{i=0}^{n-1}V(T_{X_i}\w,Z_{i+1}) - n\Lambda_a(V)} \bigr]
\]
holds and $(W_n)_{n\ge1}$ is a nonnegative martingale w.r.t. the
filtration $(\kS_0^n)_{n\ge1}$. Therefore,
%
\begin{equation}
\label{wbush} W_\infty:= \lim_{n\to\infty}W_n
\qquad\bP\mbox{-a.s.}
\end{equation}
exists. 
Moreover, the event $\{W_\infty=0\}$ is measurable w.r.t. the tail
$\s$-algebra
\[
\bigcap_{n\ge1}\s\bigl(\w_x: x\in
\bZ_+^d, \llvert x \rrvert _1 \ge n\bigr)
\]
and the Kolmogorov zero-one law implies the following dichotomy:
\begin{eqnarray*}
\mbox{either}\qquad\bP(W_\infty= 0) &=& 0\qquad\mbox{(the {\it weak disorder} regime)};
\\
\mbox{or}\qquad\bP(W_\infty= 0) &=& 1\qquad\mbox{(the {\it strong disorder}
regime)}.
\end{eqnarray*}
This analysis is due to Bolthausen \cite{Bol1989} in the case of
directed polymers (i.e., for potentials that do not depend on $z$) and
is easily adapted to our setting, which we leave to the reader. The
terms weak disorder and strong disorder were coined in \cite{ComShiYos2003}.

It follows from Jensen's inequality that $\Lambda_q(V) \le\Lambda
_a(V)$ always holds. This is known as the annealing bound. Observe
that, in the case of weak disorder, we have $\Lambda_a(V)<\infty$
(since otherwise $W_n = 0$) and
%
\begin{equation}
\label{esittir} 0 = \lim_{n\to\infty}\frac{1}{n}\log
W_n(\w) = \lim_{n\to\infty
}\frac{1}{n}\log
\cZ_{n,0}^\w- \Lambda_a(V) =
\Lambda_q(V) - \Lambda_a(V)
\end{equation}
for $\bP$-a.e. $\w$. Therefore,
\[
\Lambda_q(V) < \Lambda_a(V)\qquad\mbox{(the {\it very
strong disorder} regime)}
\]
is a sufficient condition for strong disorder. However, it is not known
whether it is necessary for strong disorder. We will say more about
this and related open problems in Section~\ref{openprobsec}.

The following theorem collects the results regarding the dependence of
the disorder regimes of directed polymers on (i) the dimension $d$ and
(ii) an inverse temperature parameter $\beta$ which is introduced to
modify the strength of the potential.

\begin{theorem}\label{critthm}
Assume that \textup{(Dir)}, \textup{(Ind)} and \textup{(Loc)} are satisfied, $V$ does not depend
on $z$, and $\Lambda_a(\beta V)<\infty$ for every $\beta\ge0$.
Then, we have the following results.
\begin{enumerate}
\item[(a)] There exist $0\le\beta_c = \beta_c(V,d) \le\beta'_c =
\beta'_c(V,d) \le\infty$ such that the RWRP (or directed polymer)
with potential $\beta V$ is in:
\begin{longlist}
\item[(i)] the weak disorder regime if $\beta\in\{0\}\cup(0,\beta_c)$,
\item[(ii)] the strong disorder regime if $\beta\in(\beta_c,\infty
)$, and
\item[(iii)] the very strong disorder regime if $\beta\in(\beta
'_c,\infty)$.
\end{longlist}
\item[(b)] The critical inverse temperatures $\beta_c = \beta
_c(V,d)$ and $\beta'_c = \beta'_c(V,d)$ satisfy
\begin{longlist}[(ii)]
\item[(i)] $\beta_c > 0$ if $d\ge4$, and
\item[(ii)] $\beta_c' = 0$ if $d=2,3$.
\end{longlist}
\end{enumerate}
\end{theorem}

Part (a) of this theorem is proved in \cite{ComYos2006}, Theorem 3.2;
item (i) of part (b) is established in a series of papers \cite
{ImbSpe1988,Bol1989,SonZho1996}; and item (ii) of part (b) is shown in
\cite{ComVar2006} for $d=2$ and \cite{Lac2010} for $d=3$. In fact,
\cite{Lac2010} covers $d=2,3$ and is valid under the weaker assumption
of $\Lambda_a(\beta_oV) < \infty$ for some $\beta_o>0$. As far as
we know, these results have not been adapted to the RWRP model with
potentials that depend on $z$. However, the analogs of items (i) and
(ii) of part (b) have been established in \cite{Yil2009a,YilZei2010}
in the context of large deviations for directed RWRE.

%

\subsection{Organization of the article}

In Section~\ref{resultssec}, we present our results along with remarks
and open problems. The subsequent sections contain the proofs of our results.

%

\section{Results}\label{resultssec}

\subsection{Quenched free energy in the general case}
In order to abbreviate the variational formula (\ref{bboduk}) given in
Theorem~\ref{KVarOld} for the quenched free energy $\Lambda_q(V)$, we define
\[
K(V,F):= \bP\mbox{-}\esssup_{\w} \biggl\{\log \biggl(\sum
_{z\in
\cR}p(z)e^{V(\w,z) + F(\w,z)} \biggr) \biggr\}
\]
for every measurable function $F:\Omega\times\cR\to\bR$. Observe
that $K(V,F)$ is equal to
\[
K'(V,g):= \bP\mbox{-}\esssup_{\w} \biggl\{\log \biggl(
\sum_{z\in
\cR}\frac{p(z)e^{V(\w,z)}g(T_z\w)}{g(\w)} \biggr) \biggr\}
\]
when $F$ is of the form
%
\begin{equation}
\label{elvsok} F(\w,z) = \bigl(\nabla^*g\bigr) (\w,z):= \log \biggl(
\frac{g(T_z\w)}{g(\w
)} \biggr)
\end{equation}
for some $g\in L^+(\Omega,\kS,\bP)$. Here and throughout,
\begin{eqnarray*}
L^+\bigl(\Omega,\kS',\bP\bigr) &:=& \bigl\{g:\Omega\to\bR:
\mbox{$g$ is $\kS '$-measurable and } 0<g(\w)<\infty \mbox{ for }\bP\mbox{-a.e. }\w\bigr\} \quad\mbox{and}
\\
L^{++}\bigl(\Omega,\kS',\bP\bigr) &:=& \bigl\{g:\Omega\to
\bR: \mbox{$g$ is $\kS '$-measurable and } \exists c>0 \mbox{ s.t. }
c<g(\w)<\infty \mbox{ for }\bP\mbox{-a.e. }\w\bigr\}
\end{eqnarray*}
for every $\s$-algebra $\kS'\subset\kS$ on $\Omega$. 
We start our analysis by showing that the logarithmic gradient (as in
(\ref{elvsok})) of any $g\in L^+(\Omega,\kS,\bP)$ is in $\cK_0$
whenever $K'(V,g)<\infty$, see Lemma~\ref{huyop}. Then, we give a
short alternative proof of (\ref{bboduk}) and provide two modified
versions of it. 

\begin{theorem}\label{skilic}
Assume that $V\in\cL$. Then, we have the following variational formulas.
%
{\renewcommand{\theequation}{qVar\arabic{equation}}\setcounter{equation}{-1}
\begin{eqnarray}
\Lambda_q(V) &=& \inf_{F\in\cK_0}K(V,F), \label{qVar0}
\\
\Lambda_q(V) &=& \inf_{g\in L^+}K'(V,g),\label{qVar1}
\\
\Lambda_q(V) &=& \inf_{g\in L^{++}}K'(V,g).\label{qVar2}
\end{eqnarray}}%
Here, the spaces $L^+$ and $L^{++}$ stand for \textup{(i)} $L^+(\Omega,\kS,\bP
)$ and $L^{++}(\Omega,\kS,\bP)$ in the general case and \textup{(ii)}
$L^+(\Omega,\kS_0^\infty,\bP)$ and $L^{++}(\Omega,\kS_0^\infty
,\bP)$  under \textup{(Dir)} and \textup{(Loc)}.
\end{theorem}

\begin{remark}\label{boyunonem}
It is shown in \cite{RasSepYil2013}, Lemma C.3, that $\cK_0$ is the
$L^1(\Omega,\kS,\bP)$-closure of
\[
\bigl\{\nabla^*g: \exists C>0 \mbox{ s.t. } C^{-1}<g(\w)<C \mbox{ for $\bP$-a.e. $\w$}\bigr\}.
\]
Unfortunately, our understanding of $\cK_0$ does not go much beyond
this characterization. Thus, for applications, \textup{(\ref{qVar0})} is perhaps not
very useful. \textup{(\ref{qVar1})} and \textup{(\ref{qVar2})} replace $\cK_0$ by the much more
concrete class of logarithmic gradients. This way, they simplify
\textup{(\ref{qVar0})} and thereby improve our understanding of the large deviation
rate functions $\kI_q$ and $I_q$ via (\ref{level2RF}) and (\ref
{level1RF}), respectively.
\end{remark}

The proof of Theorem~\ref{skilic} does not rely on the rather
technical minimax approach taken in \cite{RasSepYil2013}, Theorem 2.3.
The lower bounds in \textup{(\ref{qVar0})}, \textup{(\ref{qVar1})} and \textup{(\ref{qVar2})} follow from a standard
spectral argument, whereas the upper bounds hinge on a certain control
on the minima of path integrals of centered cocycles on large sets
which is implied by an ergodic theorem and is trivial in the case of
\textup{(\ref{qVar2})}. Moreover, as we have recorded in Theorem~\ref{pe2pe}, the
existence of the a.s. limit in (\ref{freirg}) is shown in \cite{RasSep2014}, Theorem 2.2(b), for $V\in\cL$ (without assuming that
$\kS$ is countably generated) by subadditivity and elementary
estimates. In short, the proof of Theorem~\ref{skilic} is completely
independent of Theorem~\ref{KVarOld}.

Now that we have three closely related variational formulas for
$\Lambda_q(V)$, it is natural to ask whether they possess minimizers,
that is, the infima in their definitions are attained. We provide a
positive answer to this question for \textup{(\ref{qVar0})}. As its proof in
Section~\ref{quesec1} attests, the technical significance of this
result is due to the lack of weak compactness of the unit ball in
$L^1(\Omega,\kS,\bP)$.

\begin{theorem}\label{cvtcvt}
Assume that $V\in\cL$. Then, \textup{(\ref{qVar0})} always has a minimizer.
\end{theorem}

It turns out that, unlike \textup{(\ref{qVar0})}, the variational formulas \textup{(\ref{qVar1})} and
\textup{(\ref{qVar2})} do not always have minimizers. In fact, this is one of the main
results in this paper, see Section~\ref{serkilannsec}. The possible
lack of minimizers might be seen as a shortcoming of our formulas.
However, we will argue that it is actually an advantage since it
carries valuable information about the disorder regime of the model, at
least in the directed i.i.d. case.

%

\subsection{Annealed free energy in the directed i.i.d. case}

Assume that \textup{(Dir)}, \textup{(Ind)} and \textup{(Loc)} hold. Our analysis of the
variational formulas \textup{(\ref{qVar1})} and \textup{(\ref{qVar2})} for $\Lambda_q(V)$ builds on
its analog for the annealed free energy $\Lambda_a(V)$ defined in
(\ref{avratdef}).

\begin{theorem}\label{prev}
Assume \textup{(Dir)}, \textup{(Ind)}, and \textup{(Loc)}. Then, we have the following variational
formulas.
%
{\renewcommand{\theequation}{aVar\arabic{equation}}\setcounter{equation}{0}\begin{eqnarray}
\Lambda_a(V) &=& \inf_{g\in L^+\cap L^1}K'(V,g),\label{aVar1}
\\
\Lambda_a(V) &=& \inf_{g\in L^{++}\cap L^1} K'(V,g)\label{aVar2}.
\end{eqnarray}}%
Here, $L^+$, $L^{++}$ and $L^1$ stand for $L^+(\Omega,\kS_0^\infty
,\bP)$, $L^{++}(\Omega,\kS_0^\infty,\bP)$ and $L^1(\Omega,\kS
_0^\infty,\bP)$, respectively.
\end{theorem}

\begin{remark}
The variational formulas \textup{(\ref{qVar1})} and \textup{(\ref{aVar1})} for $\Lambda_q(V)$ and
$\Lambda_a(V)$ can be equivalently written as the infima of $K(V,F)$ over
\[
\bigl\{F\in\cK_0: F = \nabla^*g \mbox{ for some } g\in L^+\bigr\}
\quad\mbox{and}\quad\bigl\{F\in\cK_0: F = \nabla^*g\mbox{ for some }
g\in L^+\cap L^1\bigr\},
\]
respectively. The presence of these different sets is not merely a
technical artifact of our proofs, as we know that $\Lambda_q(V) <
\Lambda_a(V)$ in the case of very strong disorder, cf. Theorem~\ref
{critthm}. We find this strict inequality to be particularly
interesting because both of these sets are dense in $\cK_0$ by \cite{RasSepYil2013}, Lemma C.3, cf. Remark~\ref{boyunonem}. The same
comment applies to \textup{(\ref{qVar2})} and \textup{(\ref{aVar2})}.
\end{remark}

In the light of Theorem~\ref{cvtcvt} and the paragraph below it, we
ask if/when \textup{(\ref{aVar1})} and \textup{(\ref{aVar2})} have any minimizers. The answer to this
question constitutes our first variational result on the disorder
regimes of RWRP.

\begin{theorem}\label{unique1}
Assume \textup{(Dir)}, \textup{(Ind)}, \textup{(Loc)}, and $\Lambda_a(V)<\infty$.
\begin{longlist}[(a)]
\item[(a)] \textup{(\ref{aVar1})} has a minimizer if and only if there is weak
disorder. In this case, the minimizer is unique (up to a multiplicative
constant), equal to $W_\infty$ defined in (\ref{wbush}), and there is
no need for taking essential supremum in $K'(V,W_\infty)$, that is,
\[
\Lambda_a(V) = K'(V,W_\infty) = \log \biggl(\sum
_{z\in\cR}\frac{p(z)e^{V(\w
,z)}W_\infty(T_z\w)}{W_\infty(\w)} \biggr)\qquad\mbox{for }\bP\mbox{-a.e. }\w.
\]
\item[(b)] \textup{(\ref{aVar2})} has no minimizers unless $\cZ_{1,0}^\w$ is $\bP
$-essentially constant, cf. Remark~\ref{esolmaz}.
\end{longlist}
\end{theorem}

\begin{remark}\label{esolmaz}
Theorem~\ref{unique1}(a) implies that the only minimizer candidate of
\textup{(\ref{aVar2})} is $W_\infty$. However, we will show in Proposition~\ref
{birtekbu} that $W_\infty\notin L^{++}$ unless $\cZ_{1,0}^\w= \sum_{z\in\cR}p(z)e^{V(\w,z)}$ is $\bP$-essentially constant. In the
latter case, $\cZ_{n,0}^\w$ is $\bP$-essentially\vspace*{1pt} constant for every
$n\ge1$, and $\bP(W_n = 1) = \bP(W_\infty= 1) = 1$. By Theorems
\ref{unique1}(a) and~\ref{unique2}(a), $W_\infty$ is the unique
minimizer of \textup{(\ref{aVar1})}, \textup{(\ref{aVar2})}, \textup{(\ref{qVar1})} and \textup{(\ref{qVar2})}. Observe that, in
this case, the RWRP is nothing but an RWRE with transition kernel $\hat
p(\w,z) = p(z)e^{V(\w,z) - \Lambda_a(V)}$, cf. Remark~\ref{gerektiyaw}.
\end{remark}

Other characterizations of weak disorder have been previously given in
the literature on directed polymers. First of all, it is shown in \cite{ComShiYos2003}, Theorem 2.1, that weak disorder is equivalent to the
delocalization of the polymer in an appropriate sense. Precisely, when
$\Lambda_a(V)<\infty$ and $V$ is not $\bP$-essentially constant,
there is weak disorder if and only if
\[
\sum_{n=1}^\infty\bigl(Q_{n,0}^\w
\bigr)^{\otimes2}(X_n = \tilde X_n) < \infty.
\]
Here, $\tilde X_n$ is an independent copy of $X_n$ under $Q_{n,0}^\w$.
Second, \cite{ComYos2006}, Proposition 3.1, collects some useful
characterizations of weak disorder, e.g., the $L^1(\bP)$-convergence
or uniform integrability of the martingale $(W_n)_{n\ge1}$. As far as
we know, part (a) of Theorem~\ref{unique1} is the first variational
characterization of weak disorder for RWRP. Its proof builds on an
earlier characterization given as part of \cite{ComYos2006}, Proposition~3.1, for directed polymers, see Section~\ref{issiksec} for details.

%

\subsection{Analysis of \texorpdfstring{(\protect\ref{qVar1})}{(qVar1)} and \texorpdfstring{(\protect\ref{qVar2})}{(qVar2)} in the directed i.i.d. case}\label{serkilannsec}

We continue working under \textup{(Dir)}, \textup{(Ind)} and \textup{(Loc)}. In the case of weak
disorder, $\Lambda_q(V) = \Lambda_a(V)<\infty$ by (\ref{esittir}).
Therefore, the unique minimizer $W_\infty$ of \textup{(\ref{aVar1})} in $L^+\cap L^1$
is also a minimizer of \textup{(\ref{qVar1})} in the larger space $L^+$. However, it
is not a-priori clear whether $W_\infty$ is the unique minimizer of
\textup{(\ref{qVar1})}. The following theorem settles this issue.

\begin{theorem}\label{unique2}
Assume \textup{(Dir)}, \textup{(Ind)}, \textup{(Loc)}, and weak disorder.
\begin{longlist}[(a)]
\item[(a)] Up to a multiplicative constant, the unique minimizer
$W_\infty$ of \textup{(\ref{aVar1})} is also the unique minimizer of \textup{(\ref{qVar1})}.
\item[(b)] \textup{(\ref{qVar2})} has no minimizers unless $\cZ_{1,0}^\w$ is $\bP
$-essentially constant, cf. Remark~\ref{esolmaz}.
\end{longlist}
\end{theorem}

\setcounter{equation}{1}

Note that Theorem~\ref{unique2} does not say anything about whether
\textup{(\ref{qVar1})} and \textup{(\ref{qVar2})} have any minimizers in the case of strong disorder.
This turns out to be a more difficult question. In order to address it,
we introduce
%
\begin{equation}
\label{sakinol} h_n^\lambda(\w):= E_0
\bigl[e^{\sum_{i=0}^{n-1}V(T_{X_i}\w
,Z_{i+1}) - n\lambda} \bigr]
\end{equation}
for every $n\ge1$, $\lambda\in\bR$ and $\w\in\Omega$, and
consider the future measurable functions
\[
\underline h_\infty^\lambda(\w):= \liminf_{n\to\infty
}h_n^\lambda(
\w)\quad\mbox{and}\quad\bar h_\infty^\lambda(\w):= \limsup
_{n\to\infty}h_n^\lambda(\w).
\]
With this notation, $W_n = h_n^\lambda$ and $W_\infty= \underline
h_\infty^\lambda= \bar h_\infty^\lambda$ when $\lambda= \Lambda
_a(V) < \infty$. For general $\lambda\in\bR$, we know that
\[
\lim_{n\to\infty}\frac{1}{n}\log
h_n^\lambda(\w) = \Lambda_q(V) - \lambda
\]
holds for $\bP$-a.e. $\w$. Therefore,
\begin{eqnarray*}
\bP\bigl(\underline h_\infty^\lambda= \bar
h_\infty^\lambda= 0\bigr) &=& 1\qquad\mbox{if }\lambda>\Lambda_q(V)\quad\mbox{and}
\\
\bP\bigl(\underline h_\infty^\lambda= \bar
h_\infty^\lambda= \infty\bigr) &=& 1\qquad\mbox{if }\lambda<\Lambda_q(V).
\end{eqnarray*}
Hence, the only nontrivial choice of parameter is $\lambda= \Lambda
_q(V)$. In the latter case, each of the events
%
\begin{eqnarray}
\label{olaylar}
&& \bigl\{\underline h_\infty^\lambda= 0\bigr\},\qquad
\bigl\{0 < \underline h_\infty ^\lambda< \infty\bigr\},\qquad
\bigl\{\underline h_\infty^\lambda= \infty\bigr\},
\nonumber\\[-8pt]\\[-8pt]\nonumber
&& \bigl\{ \bar h_\infty^\lambda= 0\bigr\},\qquad
\bigl\{0 < \bar h_\infty^\lambda< \infty\bigr\},\qquad
\bigl\{\bar h_\infty^\lambda= \infty\bigr\}
\end{eqnarray}
has $\bP$-probability zero or one, see Lemmas~\ref{Kolpa1} and~\ref
{Kolpa2}. To provide some insight, we make a slight digression from the
variational analysis and use one of these events to give a quenched
characterization of weak disorder.

\begin{theorem}\label{unique3}
Assume \textup{(Dir)}, \textup{(Ind)}, \textup{(Loc)}, $V\in\cL$, and $\Lambda_q(V)<\infty$.
Then, there is weak disorder if and only if $\bP(0 < \underline
h_\infty^\lambda< \infty) = 1$, i.e., $\underline h_\infty^\lambda
\in L^+$,
for $\lambda= \Lambda_q(V)$.
\end{theorem}

Next, we use another event in (\ref{olaylar}) to conditionally prove
that \textup{(\ref{qVar1})} and \textup{(\ref{qVar2})} do not always have any minimizers under strong
disorder.

\begin{theorem}\label{nominimizer2}
Assume \textup{(Dir)}, \textup{(Ind)}, \textup{(Loc)}, $V\in\cL$, and $\Lambda_q(V)<\infty$.
If there is strong disorder and $\bP(\bar h_\infty^\lambda= 0) = 0$
for $\lambda= \Lambda_q(V)$, then \textup{(\ref{qVar1})} and \textup{(\ref{qVar2})} have no minimizers.
\end{theorem}

Finally, we provide a sufficient condition for the key hypothesis of
Theorem~\ref{nominimizer2}. To this end, we fix $\lambda= \Lambda
_q(V) < \infty$ and let
%
\begin{equation}
\label{sahann} H_n(\w):= E_0 \bigl[e^{\sum_{i=0}^{n-1}V(T_{X_i}\w,Z_{i+1}) -
n\Lambda_q(V)}
\one_{\{X_n = (n/d,\ldots,n/d)\}} \bigr]
\end{equation}
be the ``bridge'' analog of $h_n^\lambda$. (For convenience, we assume
that $n$ is divisible by $d$.) With this notation, we clearly have
$h_n^\lambda\ge H_n$.

\begin{proposition}\label{teo}
Assume \textup{(Dir)}, \textup{(Ind)}, \textup{(Loc)}, $V\in\cL$, and $\Lambda_q(V)<\infty$.
If there exists an increasing sequence $(a(n))_{n\ge1}$ such that
%
\begin{equation}
\label{kondisyon} \lim_{n\to\infty}a(n) = \infty,\qquad\lim
_{n\to\infty}\frac
{a(n-1)}{a(n)} = 1\quad\mbox{and}\quad\limsup
_{n\to\infty}\bP \bigl(\log H_n \ge a(n) \bigr)>0,
\end{equation}
then
%
\begin{equation}
\label{borneo} \bP \biggl(\limsup_{n\to\infty}\frac{\log H_n}{a(n)} \ge1
\biggr) = 1.
\end{equation}
In particular, $\bP(\bar h_\infty^\lambda= \infty) =1$ for $\lambda
= \Lambda_q(V)$.
\end{proposition}

It has been recently shown in \cite{BorCorRem2013} that $n^{-1/3}\log
H_n$ has an $F_{\mathrm{GUE}}$ distributional limit for the log-gamma
directed polymer model on $\bZ^2$ with parameter $\gamma\in(0,\gamma
^*)$ for some $\gamma^*>0$. In particular, the conditions in
Proposition~\ref{teo} are satisfied with $a(n) = n^{1/3}$. On the
other hand, since $d=2$ in this example, it is in the very strong
disorder regime by \cite{ComVar2006,Lac2010}. (Technically, [10]
assumes that $\Lambda_a(\beta V)<\infty$ for every $\beta>0$, and
[25] weakens this assumption to $\Lambda_a(\beta_o V)<\infty$ for
some $\beta_o>0$. The log-gamma model satisfies only this weaker
condition.) We thereby conclude that \textup{(\ref{qVar1})} and \textup{(\ref{qVar2})} do not always
have any minimizers in the case of very strong disorder. We record this
as a remark for future reference.

\begin{remark}\label{kolonyagi}
Assume \textup{(Dir)}, \textup{(Ind)}, \textup{(Loc)}, $V\in\cL$, and $\Lambda_q(V)<\infty$.
Then, as explained in the paragraph above, \textup{(\ref{qVar1})} and \textup{(\ref{qVar2})} do not
always have any minimizers in the case of very strong disorder.
\end{remark}

%

\subsection{Additional remarks and open problems}\label{openprobsec}

We know from Theorem~\ref{critthm} that the critical inverse
temperatures $\beta_c = \beta_c(V,d)$ and $\beta'_c = \beta
'_c(V,d)$ satisfy $\beta_c = \beta'_c = 0$ for $d=2,3$, and it is
natural to expect that $\beta_c = \beta'_c$ for every $d\ge2$.
However, this is an open problem, see \cite{ComYos2006}, Remark 3.2.
Furthermore, it is generally believed that there is strong disorder at
$\beta_c$ for $d\ge4$. The latter claim is supported by the analogous
result in the context of directed polymers on trees which follows from
\cite{KahPey1976}.

With this background, here is our conjecture regarding the very strong
disorder regime and the events in (\ref{olaylar}).

\begin{conjecture}\label{kancik}
Assume \textup{(Dir)}, \textup{(Ind)}, \textup{(Loc)}, $V\in\cL$, and $\Lambda_q(V)<\infty$. Then,
\[
\bP\bigl(0 = \underline h_\infty^\lambda< \bar
h_\infty^\lambda= \infty \bigr) = 1
\]
for $\lambda= \Lambda_q(V)$ whenever there is very strong disorder.
\end{conjecture}

If this conjecture is indeed true, it would readily give the following
quenched characterization of the disorder regimes.
\begin{longlist}[(a)]
\item[(a)] If there is weak disorder, then
\[
\bP\bigl(0 < \underline h_\infty^\lambda= \bar
h_\infty^\lambda< \infty \bigr) = 1\qquad\mbox{for }\lambda=
\Lambda_q(V) = \Lambda_a(V) < \infty.
\]
\item[(b)] If there is critically strong disorder, then
\[
\bP\bigl(\underline h_\infty^\lambda= \bar h_\infty^\lambda=
0\bigr) = 1\qquad\mbox{for }\lambda= \Lambda_q(V) =
\Lambda_a(V) < \infty.
\]
\item[(c)] If there is very strong disorder, then
\[
\bP\bigl(0 = \underline h_\infty^\lambda< \bar
h_\infty^\lambda= \infty \bigr) = 1\qquad\mbox{for }\lambda=
\Lambda_q(V) < \Lambda_a(V) \le \infty.
\]
\end{longlist}
This result would constitute a stronger version of Theorem~\ref
{unique3}. Note that parts (a) and (b) are trivial since $\underline
h_\infty^\lambda= \bar h_\infty^\lambda= W_\infty$ for $\lambda=
\Lambda_a(V)$.

As a second application, if Conjecture~\ref{kancik} is true, then very
strong disorder would imply the hypotheses of Theorem~\ref
{nominimizer2}, and \textup{(\ref{qVar1})} and \textup{(\ref{qVar2})} would never have any minimizers
in that case. In other words, we could establish a stronger version of
Remark~\ref{kolonyagi}.

The result of Borodin \textit{et al}. \cite{BorCorRem2013} that we have used
to satisfy the conditions of Proposition~\ref{teo} is a form of
Kardar--Parisi--Zhang (KPZ) universality and is expected to hold for a
large class of models, see \cite{Cor2012} for a survey. However,
Proposition~\ref{teo} is much more modest since it does not require
any sharp estimates such as the $n^{1/3}$ scaling in KPZ universality.
Indeed, slowly growing sequences, for example, $a(n) = \log\log\log
n$, satisfy the first two conditions in (\ref{kondisyon}).

Finally, observe that Theorem~\ref{nominimizer2} is not applicable in
the (hypothetical) case of critically strong disorder since, then, $\bP
(\bar h_\infty^\lambda= 0) = 1$ for $\lambda= \Lambda_q(V) =
\Lambda_a(V)$. Therefore, we refrain from making any claims regarding
the existence of any minimizers of \textup{(\ref{qVar1})} and \textup{(\ref{qVar2})} in that case.

%

\section{Quenched free energy in the general case}\label{quesec}

\subsection{Variational formulas \texorpdfstring{(\protect\ref{qVar0})}{(qVar0)}, \texorpdfstring{(\protect\ref{qVar1})}{(qVar1)} and 
\texorpdfstring{(\protect\ref{qVar2})}{(qVar2)} for \texorpdfstring{$\Lambda_q(V)$}{Lambdaq(V)}}\label{yetissec}

\begin{lemma}\label{huyop}
Assume that $V(\cdot,z)\in L^1(\Omega,\kS,\bP)$ for every $z\in\cR
$. If $K(V,F)<\infty$ with $F = \nabla^*g$ as defined in (\ref
{elvsok}) for some $g\in L^+(\Omega,\kS,\bP)$, 
then $F\in\cK_0$.
\end{lemma}

\begin{pf}
It follows from the definition of $K(V,F)$ that
%
\begin{equation}
\label{daglars} F(\cdot,z) \le\bigl\llvert V(\cdot,z) \bigr\rrvert + \log\llvert
\cR \rrvert + K(V,F)
\end{equation}
$\bP$-a.s. for every $z\in\cR$. Therefore, $F^+(\cdot,z)$ is
integrable and $\bE[F(\cdot,z)]$ is well defined, even though it
might a-priori be $-\infty$. Note that $\bE[F(\cdot,z)] = -\infty$
is equivalent to $\bE[\llvert  F(\cdot,z) \rrvert  ] = \infty$.

As a consequence of telescoping, we have
%
\begin{equation}
\frac{1}{n}\sum_{i=0}^{n-1}F(T_{iz}
\w,z) = \frac{1}{n}\log\frac
{g(T_{nz}\w)}{g(\w)} = \frac{1}{n}\log
g(T_{nz}\w) - \frac{1}{n}\log g(\w).\label{emrr}
\end{equation}
By Birkhoff's ergodic theorem, the LHS of (\ref{emrr}) converges $\bP
$-a.s. (and hence also in $\bP$-pro\-bability) to $\bE[F(\cdot,z)]
\in[-\infty,\infty)$. However, the RHS of (\ref{emrr}) converges to
$0$ in $\bP$-probability. Indeed, for every $\e>0$,
\[
\bP \biggl(\biggl\llvert \frac{1}{n}\log g\circ T_{nz}\biggr
\rrvert >\e \biggr) = \bP \biggl(\biggl\llvert \frac{1}{n}\log g\biggr
\rrvert >\e \biggr) = \bP\bigl(\llvert \log g\rrvert >n\e\bigr)\to0\qquad\mbox{as }n\to\infty.
\]
We conclude that $\bE[F(\cdot,z)] = 0$ and $F(\cdot,z)\in L^1(\bP
)$. Finally, the cocycle property is obvious from the definition of
$F$. This finishes the proof.
\end{pf}

\begin{pf*}{Proof of Theorem~\ref{skilic} (The upper bounds)}
We start by considering \textup{(\ref{qVar0})}.
Take any $F\in\cK_0$ and assume WLOG that $K(V,F) < \infty$ since
the desired upper bound is otherwise trivial.
Observe that
\begin{eqnarray*}
&&E_0 \bigl[e^{\sum_{i=0}^{n-1}V(T_{X_i}\w,Z_{i+1}) + F(T_{X_i}\w
,Z_{i+1})} \bigr]
\\
&&\quad  = \sum_{x\in D_{n-1}}E_0
\bigl[e^{\sum_{i=0}^{n-2}V(T_{X_i}\w
,Z_{i+1}) + F(T_{X_i}\w,Z_{i+1})}\one_{\{X_{n-1} = x\}} \bigr] \sum
_{z\in\cR}p(z)e^{V(T_x\w,z) + F(T_x\w,z)}
\\
&&\quad  \le\sum_{x\in D_{n-1}}E_0
\bigl[e^{\sum_{i=0}^{n-2}V(T_{X_i}\w
,Z_{i+1}) + F(T_{X_i}\w,Z_{i+1})}\one_{\{X_{n-1} = x\}} \bigr]e^{K(V,F)}
\\
&&\quad  = E_0 \bigl[e^{\sum_{i=0}^{n-2}V(T_{X_i}\w,Z_{i+1}) + F(T_{X_i}\w
,Z_{i+1})} \bigr]e^{K(V,F)}\le\cdots\le
e^{nK(V,F)},
\end{eqnarray*}
where $D_{n-1}$ is defined in (\ref{mrdj}). Therefore,
%
\begin{eqnarray}
\label{yetdost} \Lambda_q(V) &=& \lim_{n\to\infty}
\frac{1}{n}\log\cZ_{n,0}^\w\nn
\nonumber\\[-8pt]\\[-8pt]\nonumber
&\le&\limsup_{n\to\infty}\frac{1}{n}\log E_0
\bigl[e^{\sum
_{i=0}^{n-1}V(T_{X_i}\w,Z_{i+1}) + F(T_{X_i}\w,Z_{i+1})} \bigr] \le K(V,F)
\end{eqnarray}
if
%
\begin{equation}
\label{ifoz} \liminf_{n\to\infty}\min_{x\in D_n}
\frac{1}{n}\sum_{i=0}^{n-1}F(T_{x_i}
\w,z_{i+1}) \ge0.
\end{equation}
Here, by the cocycle property, $(x_i)_{i=0}^n$ is allowed to be any
path such that $z_{i+1} = x_{i+1} - x_i\in\cR$, $x_0 = 0$ and $x_n =
x$. 

We see from (\ref{daglars}) that $F$ is $\bP$-a.s. bounded from
above by a function in $\cL$.
Under this assumption, it has been recently shown in \cite{GeoRasSepYil2014}, Theorem~9.3, that
\[
\lim_{n\to\infty}\max_{x\in D_n}\frac{1}{n}\Biggl
\llvert \sum_{i=0}^{n-1}F(T_{x_i}
\w,z_{i+1})\Biggr\rrvert = 0.
\]
%
This is an ergodic theorem for cocycles. In particular, we have (\ref
{ifoz}), and therefore, (\ref{yetdost}). Taking infimum over all $F\in
\cK_0$ gives the upper bound in \textup{(\ref{qVar0})}.

The upper bounds in \textup{(\ref{qVar1})} and \textup{(\ref{qVar2})} are now easy. Indeed, take any
$g\in L^+(\Omega,\kS,\bP)$ and assume WLOG that $K'(V,g) < \infty$
since the desired upper bounds are otherwise trivial.
Then, $F:= \nabla^*g\in\cK_0$ by Lemma~\ref{huyop}, and $\Lambda
_q(V) \le K(V,F) = K'(V,g)$ by the upper bound in \textup{(\ref{qVar0})}. Taking
infimum over all $g\in L^+(\Omega,\kS,\bP)$ gives the upper bound in
\textup{(\ref{qVar1})}, from which the upper bound in \textup{(\ref{qVar2})} follows.
\end{pf*}

\begin{remark}
Note that, for the logarithmic gradient $F = \nabla^*g$ of any $g\in
L^+(\Omega,\kS,\bP)$, the condition in (\ref{ifoz}) can be written as
%
\begin{equation}
\label{rukadam} \liminf_{n\to\infty}\min_{x\in D_n}
\frac{1}{n}\log\frac{g(T_x\w
)}{g(\w)} \ge0.
\end{equation}
Consequently, the upper bound in \textup{(\ref{qVar2})} does not rely on the
aforementioned ergodic theorem for cocycles because (\ref{rukadam}) is
obvious for $g\in L^{++}(\Omega,\kS,\bP)$.
\end{remark}

\begin{pf*}{Proof of Theorem~\ref{skilic} (The lower bounds)}
Assume WLOG that $\Lambda_q(V)<\infty$ since the desired lower bounds
are otherwise trivial. For any $\lambda> \Lambda_q(V)$ and $n\ge1$,
recall the function $h_n^\lambda$ which was introduced in (\ref
{sakinol}). Set $h_0^\lambda= 1$ as a convention and define
%
\begin{equation}
\label{ge} g_\lambda:= \sum_{n=0}^\infty
h_n^\lambda\ge1.
\end{equation}
Since
\[
\lim_{n\to\infty}\frac{1}{n}\log h_n^\lambda(
\w) = \Lambda_q(V) - \lambda<0
\]
for $\bP$-a.e. $\w$, we have $g_\lambda\in L^{++}(\Omega,\kS,\bP
)$. Moreover, under \textup{(Dir)} and \textup{(Loc)}, $g_\lambda$ is
future measurable.

Decompose $g_\lambda$ in the following way: for $\bP$-a.e. $\w$,
\begin{eqnarray*}
g_\lambda(\w) & =& 1 + \sum_{n=1}^\infty
h_n^\lambda(\w) = 1 + \sum_{n=1}^\infty
\sum_{z\in\cR}p(z)e^{V(\w,z) - \lambda
}h_{n-1}^\lambda(T_z
\w)
\\
& =& 1 + \sum_{z\in\cR}p(z)e^{V(\w,z) - \lambda}\sum
_{n=1}^\infty h_{n-1}^\lambda(T_z
\w) = 1 + \sum_{z\in\cR}p(z)e^{V(\w,z) -
\lambda}g_\lambda(T_z
\w).
\end{eqnarray*}
Rearranging this, we see that
%
\begin{equation}
\label{essgerek} \lambda= \log \biggl(\frac{e^\lambda}{g_\lambda(\w)} + \sum
_{z\in
\cR}\frac{p(z)e^{V(\w,z)}g_\lambda(T_z\w)}{g_\lambda(\w)} \biggr) > \log \biggl(\sum
_{z\in\cR}\frac{p(z)e^{V(\w,z)}g_\lambda
(T_z\w)}{g_\lambda(\w)} \biggr).
\end{equation}
Therefore, $\lambda\ge K'(V,g_\lambda)$. First, taking infimum over
all $g\in L^{++}$ and then taking infimum over $\lambda> \Lambda
_q(V)$ gives the lower bound in \textup{(\ref{qVar2})}, from which the lower bounds in
\textup{(\ref{qVar1})} and \textup{(\ref{qVar0})} follow since $\nabla^*g_\lambda\in\cK_0$ by
Lemma~\ref{huyop}.
\end{pf*}

%

\subsection{Minimizing the variational formula \texorpdfstring{(\protect\ref{qVar0})}{(qVar0)} for \texorpdfstring{$\Lambda_q(V)$}{Lambdaq(V)}}\label{quesec1}

\begin{pf*}{Proof of Theorem~\ref{cvtcvt}}
If $\Lambda_q(V) = \infty$, then every $F\in\cK_0$ is trivially a
minimizer of \textup{(\ref{qVar0})}. Therefore, in the rest of the proof, we will
assume that $\Lambda_q(V) < \infty$.

Since \textup{(\ref{qVar0})} involves an infimum, for every $i\ge1$, there exists an
$F_i\in\cK_0$ such that
\[
\sum_{z\in\cR}p(z)e^{V(\cdot,z) + F_i(\cdot,z)} \le e^{\Lambda
_q(V) + 1/i}
\]
holds $\bP$-a.s. Note that
\[
F_i(\cdot,z) \le\bigl\llvert V(\cdot,z) \bigr\rrvert + \log\llvert
\cR \rrvert + \Lambda_q(V) + 1/i
\]
for every $z\in\cR$. Since $V(\cdot,z)$ is in $L^1(\bP)$, we see
that $F^+_i(\cdot,z)$ is uniformly integrable. The $F_i$ are centered
by definition, so we have
$\bE[F^-_i(\cdot,z)] = \bE[F^+_i(\cdot,z)]$. Therefore, $\bE
[F^-_i(\cdot,z)]$ is uniformly bounded. By \cite{KosVar2008}, Lemma~4.3, we can write
\[
F^-_i(\cdot,z) = \hat F^-_i(\cdot,z) +
R_i(\cdot,z),
\]
where, up to a common subsequence, $\hat F^-_i(\cdot,z)$ is uniformly
integrable and $R_i(\cdot,z)\ge0$ converges to $0$ in $\bP
$-probability. Extracting a further subsequence, $\tilde F_i(\cdot,z)
= F^+_i(\cdot,z) - \hat F^-_i(\cdot,z)$ is weakly convergent in
$L^1(\bP)$ to some $\tilde F(\cdot,z)$, and $R_i(\cdot,z)$ converges
$\bP$-a.s. to $0$. By \cite{Rud1991}, Theorem 3.12, $\tilde F(\cdot
,z)$ is in the strong $L^1(\bP)$-closure of the convex hull of $\{
\tilde F_i(\cdot,z): i\ge1\}$, that is, there exists a finite convex
combination $\tilde G_i(\cdot,z):= \sum_{j=i}^\infty\alpha
_{i,j}\tilde F_j(\cdot,z)$ that converges to $\tilde F(\cdot,z)$
strongly in $L^1(\bP)$. Up to a further subsequence, $\tilde G_i(\cdot
,z)$ converges $\bP$-a.s. to $\tilde F(\cdot,z)$. This ensures that
$\tilde F(\cdot,z)$ satisfies the cocycle property. Moreover, since
$R_i(\cdot,z)\ge0$, we have $c(z):= \bE[\tilde F(\cdot,z)]\ge0$.
Let $F(\cdot,z) = \tilde F(\cdot,z) - c(z)$ for every $z\in\cR$.
Then, $F\in\cK_0$. By Jensen's inequality,
\[
\sum_{z\in\cR}p(z)e^{V(\cdot,z) + \tilde G_i(\cdot,z) - \sum
_{j=i}^\infty\alpha_{i,j}R_j(\cdot,z)} \le e^{\Lambda_q(V) + 1/i}.
\]
Sending $i\to\infty$, we get
\[
\sum_{z\in\cR}p(z)e^{V(\w,z) + F(\w,z) + c(z)} \le e^{\Lambda_q(V)}
\]
for $\bP$-a.e. $\w$ and conclude that $F$ is a minimizer of \textup{(\ref{qVar0})}.
Plus, we deduce that $c(z)=0$ for every $z\in\cR$ since, otherwise,
the RHS of \textup{(\ref{qVar0})} would be strictly less than $\Lambda_q(V)$.
\end{pf*}

%

\section{Annealed free energy in the directed i.i.d. case}

In the rest of the paper, $L^+$, $L^{++}$ and $L^1$ stand for
$L^+(\Omega,\kS_0^\infty,\bP)$, $L^{++}(\Omega,\kS_0^\infty,\bP
)$ and $L^1(\Omega,\kS_0^\infty,\bP)$, respectively.

\subsection{Variational formulas \texorpdfstring{\textup{(\protect\ref{aVar1})}}{(aVar1)} and \texorpdfstring{\textup{(\protect\ref{aVar2})}}{(aVar2)} for \texorpdfstring{$\Lambda_a(V)$}{Lambdaa(V)}}

\begin{pf*}{Proof of Theorem~\ref{prev} (The upper bounds)}
Take any $g\in L^+\cap L^1$ and assume WLOG that $K'(V,g)<\infty$
since the desired upper bounds are otherwise trivial. Then, for $\bP
$-a.e. $\w$,
\[
K'(V,g) \ge\log \biggl(\sum_{z\in\cR}
\frac{p(z)e^{V(\w,z)}g(T_z\w
)}{g(\w)} \biggr).
\]
Rearranging this, we get
%
\begin{equation}
\label{kirgiz} g(\w) \ge\sum_{z\in\cR}p(z)e^{V(\w,z) - K'(V,g)}g(T_z
\w).
\end{equation}
For every $z\in\cR$, the random variables $V(\cdot,z)$ and $g\circ
T_z$ are independent by \textup{(Dir)}, \textup{(Ind)}, \textup{(Loc)} and the future
measurability of $g$. Taking the expectation of both sides of (\ref
{kirgiz}), we see that
%
\begin{eqnarray}\label{ulupamir}
\bE[g] &\ge&\sum_{z\in\cR}p(z)\bE \bigl[e^{V(\cdot
,z)-K'(V,g)}g
\circ T_z \bigr]\nn
\nonumber\\[-8pt]\\[-8pt]\nonumber
&=& \sum_{z\in\cR}p(z)\bE \bigl[e^{V(\cdot,z)-K'(V,g)} \bigr]
\bE [g\circ T_z ] = e^{\Lambda_a(V) - K'(V,g)}\bE[g]
\end{eqnarray}
by stationarity, which implies $\Lambda_a(V) \le K'(V,g)$. The infimum
over all $g\in L^+\cap L^1$ gives the upper bound in \textup{(\ref{aVar1})}, from
which the upper bound in \textup{(\ref{aVar2})} follows.
\end{pf*}

\begin{pf*}{Proof of Theorem~\ref{prev} (The lower bounds)}
Assume WLOG that $\Lambda_a(V)<\infty$ since the desired lower bounds
are otherwise trivial. Take any $\lambda> \Lambda_a(V)$ and recall
the function $g_\lambda\in L^{++}$ which is defined in (\ref{ge}).
Its expected value is easy to compute:
\[
\bE[g_\lambda] = \sum_{n=0}^\infty\bE
\bigl[h_n^\lambda\bigr] = \sum_{n=0}^\infty
e^{n(\Lambda_a(V)-\lambda)} = \frac{1}{1 - e^{\Lambda
_a(V)-\lambda}}<\infty.
\]
Therefore, $g_\lambda\in L^{++}\cap L^1$. We have seen in (\ref
{essgerek}) that $\lambda\ge K'(V,g_\lambda)$. Taking first  infimum
over all $g\in L^{++}\cap L^1$ and  then infimum over $\lambda>
\Lambda_a(V)$ gives the lower bound in \textup{(\ref{aVar2})}, from which the lower
bound in \textup{(\ref{aVar1})} follows.
\end{pf*}

%

\subsection{An annealed variational characterization of weak disorder}\label{issiksec}




\begin{lemma}\label{alfann}
Assume \textup{(Dir)}, \textup{(Ind)}, and \textup{(Loc)}. Then, weak disorder is equivalent to
the existence of a function $g\in L^+\cap L^1$ such that
%
\begin{equation}
\label{prop1} g = \sum_{z\in\cR}p(z)e^{V(\cdot,z)-\lambda}g
\circ T_z
\end{equation}
$\bP$-a.s. for some $\lambda\in\bR$. In that case, $\lambda=
\Lambda_a(V)$, and $g$ is equal (up to a multiplicative constant) to
$W_\infty$ which is defined in (\ref{wbush}).
\end{lemma}

\begin{remark}
This result has been previously obtained as part of \cite{ComYos2006}, Proposition~3.1, in the case of directed polymers, that is, for
potentials that do not depend on $z$. Our proof below is a
straightforward adaptation, which we include for the sake of
completeness as well as for demonstrating a technique that we will use
in the rest of the paper.
\end{remark}

\begin{pf*}{Proof of Lemma \ref{alfann}}
If there is weak disorder, then $\Lambda_a(V) < \infty$ and $W_\infty
\in L^+$ by definition. Observe that
\[
\bE[W_\infty] \le\liminf_{n\to\infty}\bE[W_n]=1
\]
by Fatou's lemma, so in fact $W_\infty\in L^+\cap L^1$.
Decompose $W_n$ with respect to the first step of the underlying random
walk and see that
%
\begin{equation}
W_n = \sum_{z\in\cR}p(z)e^{V(\cdot,z)-\Lambda_a(V)}W_{n-1}
\circ T_z.\label{useit}
\end{equation}
Taking $n\to\infty$ gives (\ref{prop1}) with $g = W_\infty$ and
$\lambda= \Lambda_a(V)$.

Conversely, if there exists some $g\in L^+\cap L^1$ and $\lambda\in
\bR$ such that (\ref{prop1}) is satisfied, then we take the
expectation of both sides of (\ref{prop1}) and get
\[
\bE[g] 
= \sum_{z\in\cR}p(z)\bE
\bigl[e^{V(\cdot,z)-\lambda} \bigr]\bE [g\circ T_z ] = e^{\Lambda_a(V) - \lambda}
\bE[g]
\]
%
which implies that $\Lambda_a(V) = \lambda< \infty$. Here, as in
(\ref{ulupamir}), we used \textup{(Dir)}, \textup{(Ind)}, \textup{(Loc)} and the future
measurability of $g$. Iterating (\ref{prop1}) for $n\ge1$ times, we get
\begin{eqnarray*}
g &=& E_0 \Biggl[\exp \Biggl(\sum_{i=0}^{n-1}V(T_{X_i}
\cdot, Z_{i+1})-n\Lambda_a(V) \Biggr)g\circ
T_{X_n} \Biggr]
\\
&=& \sum_xh_n^\lambda(
\cdot,x)g\circ T_x
\end{eqnarray*}
with $\lambda= \Lambda_a(V)$ and
\[
h_n^\lambda(\cdot,x) = E_0 \bigl[e^{\sum_{i=0}^{n-1}V(T_{X_i}\cdot,
Z_{i+1}) - n\lambda}
\one_{\{X_n = x\}} \bigr].
\]
Observe that $h_n^\lambda(\cdot,x)$ is $\kS_0^n$-measurable since
$X_n$ only takes values $x\in\bZ_+^d$ such that $\llvert  x \rrvert  _1 = n$. On the
other hand, $g\circ T_x$ is independent of $\kS_0^n$ since $g$ is
future measurable. Therefore,
%
\begin{equation}
\label{peeloff} \bE\bigl[g \mid \kS_0^n\bigr] = \sum
_xh_n^\lambda(\cdot,x)\bE[g
\circ T_x] = \sum_xh_n^\lambda(
\cdot,x)\bE[g] = W_n\bE[g].
\end{equation}
Finally,
\[
W_\infty= \lim_{n\to\infty}W_n = \lim
_{n\to\infty}\frac{\bE[g \mid  \kS_0^n]}{\bE[g]} = \frac{g}{\bE[g]}>0
\]
holds $\bP$-a.s. and we conclude that there is weak disorder.
\end{pf*}

\begin{pf*}{Proof of Theorem~\ref{unique1}}
If there is weak disorder, then by Lemma~\ref{alfann}, there exists a
$g\in L^+\cap L^1$ that satisfies (\ref{prop1}) with $\lambda=
\Lambda_a(V) < \infty$. Rearranging this equality, we immediately see
that $g$ is a minimizer of \textup{(\ref{aVar1})} and there is no need for taking
essential supremum in $K'(V,g)$.

Conversely, if $\Lambda_a(V) < \infty$ and \textup{(\ref{aVar1})} has a minimizer
$g\in L^+\cap L^1$, then we have
%
\begin{equation}
\label{noesssup} g(\w) \ge\sum_{z\in\cR}p(z)e^{V(\w,z)-\Lambda_a(V)}g(T_z
\w)
\end{equation}
for $\bP$-a.e. $\w$. If taking essential supremum in $K'(V,g)$ were
indeed necessary, then the inequality in (\ref{noesssup}) would be
strict on a set of positive $\bP$-probability. In that case, we would have
\begin{eqnarray*}
\bE[g] & > &\sum_{z\in\cR}p(z)\bE \bigl[e^{V(\cdot,z)-\Lambda
_a(V)}g
\circ T_z \bigr]
\\
& =& \sum_{z\in\cR}p(z)\bE \bigl[e^{V(\cdot,z)-\Lambda_a(V)} \bigr]
\bE [g\circ T_z ] = e^{\Lambda_a(V) - \Lambda_a(V)}\bE [g] = \bE[g]
\end{eqnarray*}
which is a contradiction. Hence, there is no need for taking essential
supremum in $K'(V,g)$. Therefore, $g$ satisfies (\ref{prop1}) with
$\lambda= \Lambda_a(V)$. By Lemma~\ref{alfann}, we have weak
disorder and $g$ is equal (up to a multiplicative constant) to
$W_\infty$. This concludes the proof of part (a).

For part (b), note that any minimizer of \textup{(\ref{aVar2})} would be a minimizer
of \textup{(\ref{aVar1})}. Therefore, by part (a), \textup{(\ref{aVar2})} has no minimizers under
strong disorder, and has at most one minimizer under weak disorder,
namely $W_\infty$. However, the latter is ruled out by Proposition
\ref{birtekbu} below unless $\cZ_{1,0}^\w$ is $\bP$-essentially constant.
\end{pf*}

\begin{proposition}\label{birtekbu}
Assume \textup{(Dir)}, \textup{(Ind)}, \textup{(Loc)}, and weak disorder. Then,
\[
\textup{(a)}\quad\bP\mbox{-}\essinf_\w W_\infty(\w) = 0\quad\mbox{and}\quad\textup{(b)}\quad \bP\mbox{-}\esssup_\w W_\infty(\w) = \infty
\]
unless
\[
\cZ_{1,0}^\w= \sum_{z\in\cR}p(z)e^{V(\w,z)}
\]
is $\bP$-essentially constant, cf. Remark~\ref{esolmaz}. 
\end{proposition}

\begin{pf}
Let us prove part (a) by contradiction. Suppose $\exists c\in(0,1)$
such that $\bP(W_\infty> c) = 1$. Then, we have $\bP(W_n \le c) = 0$
for every $n\ge1$ because, otherwise,
\begin{eqnarray*}
c\bP(W_n\le c) &<& \bE[W_\infty\one_{\{W_n\le c\}}] = \bE
\bigl[\bE \bigl[W_\infty\one_{\{W_n\le c\}} \mid \kS_0^n
\bigr]\bigr] = \bE\bigl[\bE\bigl[W_\infty \mid \kS_0^n
\bigr]\one_{\{W_n\le c\}}\bigr]
\\
&=& \bE[W_n \one_{\{W_n\le c\}}] \le c\bP(W_n\le c).
\end{eqnarray*}
On the other hand, if $\cZ_{1,0}^\w$ is not $\bP$-essentially
constant, then $\bP(W_n \le c) > 0$ for large $n\ge1$. Indeed, $\bE
[\cZ_{1,0}^\w] = e^{\Lambda_a(V)}$ and there\vspace*{1pt} exists a $\delta>0$
such that $\bP(\cZ_{1,0}^\w\le e^{\Lambda_a(V) - \delta}) > 0$. By
the assumptions \textup{(Ind)} and \textup{(Loc)}, the event
\[
\bigcap_{\llvert  x \rrvert  _1\le n-1} \bigl\{\w: \cZ_{1,x}^\w
\le e^{\Lambda_a(V)
- \delta} \bigr\}
\]
has positive $\bP$-probability. On this event,
\begin{eqnarray*}
W_n(\w) &=& \sum_x W_{n-1}(
\w,x)\sum_{z\in\cR}p(z)e^{V(T_x\w,z) -
\Lambda_a(V)} = \sum
_x W_{n-1}(\w,x)\cZ_{1,x}^\w
e^{-\Lambda
_a(V)}
\\
&\le& W_{n-1}(\w)e^{-\delta} \le\cdots\le e^{-n\delta} \le c
\end{eqnarray*}
for $n\ge{\llvert  \log c\rrvert  }/{\delta}$. Here, $W_{n-1}(\w,x):=
h_{n-1}^\lambda(\w,x)$ with $\lambda= \Lambda_a(V)$. The proof of
part (b) is similar.
\end{pf}

%

\section{Analysis of \texorpdfstring{(\protect\ref{qVar1})}{(qVar1)} and \texorpdfstring{(\protect\ref{qVar2})}{(qVar2)} in the directed i.i.d. case}

\subsection{Quenched variational analysis of weak disorder}

\begin{pf*}{Proof of Theorem~\ref{unique2}}
First, without assuming weak disorder, suppose $V\in\cL$, $\Lambda
_q(V) < \infty$, and $g\in L^+$ is a minimizer of \textup{(\ref{qVar1})}. Then, it satisfies
%
\begin{equation}
\label{kullanis} g(\w) \ge\sum_{z\in\cR}p(z)e^{V(\w,z)-\Lambda_q(V)}g(T_z
\w)
\end{equation}
for $\bP$-a.e. $\w$. Iterating this inequality for $n\ge1$ times,
we see that
\[
g(\w) \ge\sum_xh_n^\lambda(
\w,x)g(T_x\w)
\]
holds with $\lambda= \Lambda_q(V)$. Dividing both sides by
$h_n^\lambda(\w)$, we get
\[
\frac{g(\w)}{h_n^\lambda(\w)} \ge\sum_x\mu_n(
\w,x)g(T_x\w),
\]
where
\[
\mu_n(\w,x):= \frac{h_n^\lambda(\w,x)}{h_n^\lambda(\w)} = Q_{n,0}^\w(X_n
= x)
\]
does not depend on $\lambda$. For any $0<M<\infty$,
\[
\frac{g(\w)}{h_n^\lambda(\w)}\wedge M \ge\sum_x
\mu_n(\w,x) \bigl(g(T_x\w)\wedge M \bigr)
\]
by Jensen's inequality since $u\mapsto u\wedge M$ is a concave
function. Note that, as in the proof of Lemma~\ref{alfann}, for every
$x\in\bZ_+^d$ with $\llvert  x \rrvert  _1 = n$, the random variable $\mu_n(\cdot
,x)$ (resp., $g\circ T_x$) is measurable w.r.t. (resp., independent
of) the $\s$-algebra $\kS_0^n$. Therefore,
%
\begin{eqnarray}\label{carries}
\bE \biggl[\frac{g}{h_n^\lambda}\wedge M \Big| \kS _0^n
\biggr] &\ge&\sum_x\mu_n(\cdot,x)\bE
\bigl[(g\circ T_x)\wedge M\bigr]
\nonumber\\[-8pt]\\[-8pt]\nonumber
&=& \sum_x\mu_n(\cdot,x)\bE[g\wedge
M] = \bE[g\wedge M].\nn
\end{eqnarray}

In the case of weak disorder, we know that $\lambda= \Lambda_q(V) =
\Lambda_a(V)$ by (\ref{esittir}), $h_n^\lambda$ converges $\bP
$-a.s. to $W_\infty\in L^+\cap L^1$ as $n\to\infty$, and $W_\infty
$ is a minimizer of \textup{(\ref{qVar1})}. By the dominated convergence theorem for
conditional expectations (see \cite{Dur2010}, Theorem 5.5.9), the LHS
of (\ref{carries}) converges $\bP$-a.s. to $(g/W_\infty)\wedge M$
as $n\to\infty$. Therefore,
\[
\bE[g\wedge M] \le\frac{g}{W_\infty}\wedge M \le\frac{g}{W_\infty
} <\infty
\]
holds $\bP$-a.s. Sending $M\to\infty$ and applying the monotone
convergence theorem, we see that $\bE[g]<\infty$. So, $g\in L^+\cap
L^1$ and it is a minimizer of \textup{(\ref{aVar1})}. By Theorem~\ref{unique1}, $g$
is equal (up to a multiplicative constant) to $W_\infty$. This
concludes the proof of part (a). Finally, part (b) follows from
Proposition~\ref{birtekbu} since $W_\infty\notin L^{++}$ unless $\cZ
_{1,0}^\w$ is $\bP$-essentially constant.
\end{pf*}

%

\subsection{A quenched characterization of weak disorder}

\begin{lemma}\label{Kolpa1}
For $\lambda= \Lambda_q(V) < \infty$, each of the events $\{
\underline h_\infty^\lambda=0\}$, $\{0 < \underline h_\infty^\lambda
< \infty\}$ and $\{\underline h_\infty^\lambda=\infty\}$ has $\bP
$-probability zero or one.
\end{lemma}

\begin{pf}
For every $m,n\ge1$ and $x\in\bZ_+^d$ such that $\llvert  x \rrvert  _1 = m$, we have
%
\begin{equation}
h_{m+n}^\lambda= \sum_y
h_m^\lambda(\cdot,y)h_n^\lambda\circ
T_y \ge h_m^\lambda(\cdot,x)h_n^\lambda
\circ T_x,\label{basla}
\end{equation}
where the equality follows from decomposing the LHS w.r.t. the
possible values of $X_m$. Taking liminf of both sides as $n\to\infty
$, we get
\[
\underline h_\infty^\lambda\ge h_m^\lambda(
\cdot,x)\underline h_\infty^\lambda\circ T_x.
\]
Therefore,
%
\begin{equation}
\bigl\{\w: \underline h_\infty^\lambda(\w)=0\bigr\} \subset
\bigcap_{m=1}^\infty\bigcap
_{\llvert  x \rrvert  _1 = m}\bigl\{\w: \underline h_\infty^\lambda
(T_x\w)=0\bigr\}\label{timo1}
\end{equation}
and
%
\begin{equation}
\bigl\{\w: \underline h_\infty^\lambda(\w)=\infty\bigr\} \supset
\bigcup_{m=1}^\infty\bigcup
_{\llvert  x \rrvert  _1 = m}\bigl\{\w: \underline h_\infty^\lambda
(T_x\w)=\infty\bigr\}.\label{timo2}
\end{equation}
If $\bP(\underline h_\infty^\lambda= 0) <1$, then by ergodicity the
RHS of (\ref{timo1}) is a $\bP$-probability zero event and,
therefore, we in fact have $\bP(\underline h_\infty^\lambda= 0) =
0$. Similarly, if $\bP(\underline h_\infty^\lambda= \infty) > 0$,
then by ergodicity the RHS of (\ref{timo2}) is a $\bP$-probability
one event and, therefore, we in fact have $\bP(\underline h_\infty
^\lambda= \infty) = 1$.
\end{pf}

\begin{pf*}{Proof of Theorem~\ref{unique3}}
One direction is immediate. Indeed, if there is weak disorder, then
\[
\underline h_\infty^\lambda= \liminf_{n\to\infty}
h_n^\lambda= \lim_{n\to\infty} W_n =
W_\infty\in L^+
\]
for $\lambda= \Lambda_q(V) = \Lambda_a(V)$.

Conversely, assume that $\underline h_\infty^\lambda\in L^+$ for
$\lambda= \Lambda_q(V)<\infty$. Then, for every $n\ge1$, we have
\[
h_{n+1}^\lambda= \sum_{z\in\cR}p(z)e^{V(\cdot,z)-\Lambda
_q(V)}h_n^\lambda
\circ T_z.
\]
Taking liminf of both sides as $n\to\infty$, we get
\[
\underline h_\infty^\lambda\ge\sum_{z\in\cR}p(z)e^{V(\cdot
,z)-\Lambda_q(V)}
\underline h_\infty^\lambda\circ T_z.
\]
Multiplying both sides of this inequality by ${e^{\Lambda
_q(V)}}/{\underline h_\infty^\lambda}$ and then taking logarithm, we
see that $\Lambda_q(V)\ge K'(V,\underline h_\infty^\lambda)$, so
$\underline h_\infty^\lambda$ is a minimizer of \textup{(\ref{qVar1})}. The proof of
Theorem~\ref{unique2} carries over until (\ref{carries}) and we have
\[
\bE \biggl[\frac{\underline h_\infty^\lambda}{h_n^\lambda
}\wedge M \Big| \kS_0^n
\biggr] \ge\bE\bigl[\underline h_\infty ^\lambda\wedge M\bigr]
\]
for every $0<M<\infty$. Observe that
\[
\limsup_{n\to\infty}\frac{\underline h_\infty^\lambda
}{h_n^\lambda} = \biggl(\liminf
_{n\to\infty}\frac{h_n^\lambda
}{\underline h_\infty^\lambda} \biggr)^{-1} = 1.
\]
By a simple modification of the dominated convergence theorem for
conditional expectations (see Lemma~\ref{ModDurrett} below), we have
\[
\bE\bigl[\underline h_\infty^\lambda\wedge M\bigr] \le\limsup
_{n\to\infty
}\bE \biggl[\frac{\underline h_\infty^\lambda}{h_n^\lambda
}\wedge M \Big|
\kS_0^n \biggr]\le1\wedge M.
\]
We send $M\to\infty$ and get $\bE[\underline h_\infty^\lambda]\le
1$ by the monotone convergence theorem. Therefore, $\underline h_\infty
^\lambda\in L^+\cap L^1$ and it is a minimizer of \textup{(\ref{aVar1})} since
$\Lambda_a(V) \le K'(V,\underline h_\infty^\lambda) = \Lambda_q(V)
\le\Lambda_a(V)$. In particular, $\Lambda_a(V)<\infty$. Finally, we
use Theorem~\ref{unique1} to conclude that there is weak disorder.
\end{pf*}

\begin{lemma}\label{ModDurrett}
Let $Y_n, Y$ and $Z$ be future measurable functions such that
$ Y = \limsup_{n\to\infty} Y_n$, $\llvert  Y_n\rrvert  \le Z$ for all
$n\ge1$, and $\bE[Z] <\infty$. Then, $\limsup_{n\to
\infty}\bE[Y_n \mid \kS_0^n]\le Y$ holds $\bP$-a.s.
\end{lemma}

\begin{pf}
Let $U_N = \sup\{Y_n - Y: n\ge N\}$ for every $N\ge1$. Then,
$\llvert  U_N\rrvert  \le2Z$, so $\bE[\llvert  U_N\rrvert  ]<\infty$. Now,
\[
\limsup_{n\to\infty}\bE\bigl[Y_n - Y \mid
\kS_0^n\bigr] \le\lim_{n\to
\infty}\bE
\bigl[U_N \mid \kS_0^n\bigr] =
U_N.
\]
Sending $N\to\infty$, we see that
\[
\limsup_{n\to\infty}\bE\bigl[Y_n - Y \mid
\kS_0^n\bigr]\le\lim_{N\to\infty
}U_N
= \limsup_{n\to\infty}Y_n - Y = 0.
\]
We conclude that
\[
\limsup_{n\to\infty}\bE\bigl[Y_n \mid
\kS_0^n\bigr]\le\lim_{n\to\infty
}\bE\bigl[Y
\mid \kS_0^n\bigr] = Y.
\]\upqed
\end{pf}

%

\subsection{Quenched variational analysis of strong disorder}\label
{strongdisorder}

\begin{lemma}\label{Kolpa2}
For $\lambda= \Lambda_q(V) < \infty$, each of the events $\{\bar
h_\infty^\lambda=0\}$, $\{0 < \bar h_\infty^\lambda< \infty\}$ and
$\{\bar h_\infty^\lambda=\infty\}$ has $\bP$-probability zero or one.
\end{lemma}

\begin{pf}
Taking limsup as $n\to\infty$ of both sides of the inequality in
(\ref{basla}), we get
\[
\bar h_\infty^\lambda\ge h_m^\lambda(\cdot,x)
\bar h_\infty^\lambda \circ T_x
\]
for every $m\ge1$ and $x\in\bZ_+^d$ such that $\llvert  x \rrvert  _1 = m$.
Therefore, the set relations (\ref{timo1}) and (\ref{timo2}) hold
with $\underline h_\infty^\lambda$ replaced by $\bar h_\infty
^\lambda$. The rest of the proof is identical to that of Lemma~\ref{Kolpa1}.
\end{pf}

\begin{pf*}{Proof of Theorem~\ref{nominimizer2}}
Fix $\lambda= \Lambda_q(V)$ and assume that $\bP(0 < \bar h_\infty
^\lambda\le\infty)=1$. Take any future measurable function $g$
satisfying $\bP(0\le g<\infty)=1$ and (\ref{kullanis}). Our strategy
will be to show that $g\in L^1$. The proof of Theorem~\ref{unique2}
carries over until (\ref{carries}) and we have
\[
\bE \biggl[\frac{g}{h_n^\lambda}\wedge M \Big| \kS _0^n
\biggr] \ge\bE[g\wedge M]
\]
for every $0<M<\infty$. Pick a sufficiently small $\delta>0$ such
that $\bP(\bar h_\infty^\lambda> \delta) >0$. Let
\[
\n_1 = \n_1(\w) = \inf\bigl\{n\ge1: h_n^\lambda(
\w)\ge\delta\bigr\}
\]
be the first time that $h_n^\lambda(\w)\ge\delta$ (if such a time
exists, otherwise it is infinite). Similarly, for every $k\ge2$, let
\[
\n_k = \n_k(\w) = \inf\bigl\{n > \n_{k-1}:
h_n^\lambda(\w)\ge\delta \bigr\}.
\]
%
Each $\n_k$ is an $\bN\cup\{\infty\}$-valued stopping time and we
can consider the $\sigma$-algebras
\[
\kS_0^{\n_k}:= \bigl\{A\in\kS_0^\infty:
A\cap\{\n_k\le n\}\in\kS _0^n\mbox{ for every } n\ge1\bigr\}.
\]
For every $k\ge1$, we have
%
\begin{eqnarray}\label{justify2}
\bE \biggl[\frac{g}{h_{\n_k}^\lambda}\wedge M \Big|\kS _0^{\n_k}
\biggr]\one_{\{\n_k < \infty\}} & =& \sum_{n=1}^\infty
\bE \biggl[\frac{g}{h_n^\lambda}\wedge M \Big| \kS_0^{\n
_k}
\biggr]\one_{\{\n_k = n\}}\nn
\\
& =& \sum_{n=1}^\infty\bE \biggl[
\frac{g}{h_n^\lambda}\wedge M \Big| \kS_0^n \biggr]
\one_{\{\n_k = n\}}
\\
& \ge&\sum_{n=1}^\infty\bE[g\wedge M]
\one_{\{\n_k = n\}} = \bE [g\wedge M]\one_{\{\n_k < \infty\}}.\nn
\end{eqnarray}
Here, (\ref{justify2}) follows from Lemma~\ref{stoppingtime} below.
Now, on the set
\[
\bigl\{\bar h_\infty^\lambda>\delta\bigr\} \subset\bigcap
_{k\ge1}\{\n _k<\infty\},
\]
we have $h_{\n_k}^\lambda\ge\delta$ for every $k\ge1$, and therefore
\[
\bE[g\wedge M]\one_{\{\bar h_\infty^\lambda> \delta\}}\le\bE \biggl[\frac{g}{\delta}\wedge M \Big|
\kS_0^{\n
_k} \biggr]\one_{\{\bar h_\infty^\lambda> \delta\}}.
\]
Since $\kS_0^{\n_k}\uparrow\kS_0^\infty$ as $k\to\infty$, we
deduce that
%
\begin{equation}
\label{timocomment} \bE[g\wedge M] \le(g/\delta)\wedge M
\end{equation}
on the set $\{\bar h_\infty^\lambda>\delta\}$. 
Next, send $M\to\infty$ and apply the monotone convergence theorem to get
$\bE[g] \le g/\delta$ on the same set. This means that $g\in L^1$. On
the other hand, if $g$ were in $L^+\cap L^1$, then it would be a
minimizer of \textup{(\ref{aVar1})} since it would satisfy $\Lambda_a(V) \le K'(V,g)
\le\Lambda_q(V) \le\Lambda_a(V)$ by rearranging (\ref{kullanis}),
and there would be weak disorder by Theorem~\ref{unique1}. Therefore,
$g\notin L^+$. We conclude that \textup{(\ref{qVar1})} has no minimizers since every
minimizer of \textup{(\ref{qVar1})} must satisfy (\ref{kullanis}) and be in $L^+$. In
particular, \textup{(\ref{qVar2})} also has no minimizers.
\end{pf*}

\begin{remark}
We can strengthen the argument at the end of the proof of Theorem~\ref
{nominimizer2} in the following way. Since $g\notin L^+$, the set $\{g
= 0\}$ has positive $\bP$-probability. We can pick a sufficiently
small $\delta>0$ such that $\{g=0\}$ and $\{\bar h_\infty^\lambda
>\delta\}$ have a nontrivial intersection. Then, $\bE[g\wedge M] = 0$
by (\ref{timocomment}) and, sending $M\to\infty$, we see that $\bE
[g] = 0$, that is, $\bP(g=0) = 1$.
\end{remark}

\begin{lemma}\label{stoppingtime}
Let $\tau$ be an $\bN\cup\{\infty\}$-valued stopping time for the
filtration $(\kS_0^n)_{n\ge1}$. Consider the $\sigma$-algebra
\[
\kS_0^{\tau}:= \bigl\{A\in\kS_0^\infty:
A\cap\{\tau\le n\}\in\kS _0^n \mbox{ for every } n\ge1
\bigr\}.
\]
Then, for every $Y\in L^1(\Omega,\kS_0^\infty,\bP)$ and $n\ge1$,
\[
\bE\bigl[Y \mid \kS_0^{\tau}\bigr]\one_{\{\tau= n\}} =
\bE\bigl[Y\one_{\{\tau=
n\}} \mid \kS_0^{\tau}\bigr] = \bE
\bigl[Y \mid \kS_0^n\bigr]\one_{\{\tau= n\}}.
\]
\end{lemma}

\begin{pf}
The first equality is automatic since $\{\tau=n\}$ is $\kS_0^\tau
$-measurable. For the second equality, start by noting that $\bE[Y \mid
\kS_0^n]\one_{\{\tau= n\}}$ is $\kS_0^{\tau}$-measurable. Indeed,
for $m\ge n$,
\[
\bE\bigl[Y \mid \kS_0^n\bigr]\one_{\{\tau= n\}}
\one_{\{\tau\le m\}} = \bE \bigl[Y \mid \kS_0^n\bigr]
\one_{\{\tau= n\}}
\]
is $\kS_0^n$- and hence $\kS_0^m$-measurable. The case $m < n$ is
trivial since then
\[
\bE\bigl[Y \mid \kS_0^n\bigr]\one_{\{\tau= n\}}
\one_{\{\tau\le m\}} = 0.
\]
Now, for every bounded $\kS_0^{\tau}$-measurable test function
$\varphi$,
\[
\bE[Y\one_{\{\tau= n\}} \varphi] = \bE\bigl[ \bE\bigl[Y \mid
\kS_0^n\bigr] \one_{\{\tau= n\}} \varphi\bigr]
\]
since $\varphi\one_{\{\tau= n\}}$ is $\kS_0^n$-measurable. This
concludes the proof.
\end{pf}

%

\subsection{A sufficient condition for the nonexistence of minimizers of \texorpdfstring{(\protect\ref{qVar1})}{(qVar1)} and \texorpdfstring{(\protect\ref{qVar2})}{(qVar2)}}\label{asympsect}

\begin{pf*}{Proof of Proposition~\ref{teo}}
It is clear from (\ref{sahann}) that
\[
\log H_{m+n} \ge\log H_m + \log H_n\circ
T_{(m/d,\ldots,m/d)}
\]
for every $m,n\ge1$ (and divisible by $d$). By \textup{(Dir)}, \textup{(Ind)} and \textup{(Loc)},
the summands on the RHS are independent. 
Let
\[
\e= \limsup_{n\to\infty}\bP\bigl(\log H_n
\ge a(n)\bigr) > 0.
\]
Fix an arbitrary $\delta\in(0,1)$. There exists an $m_1\ge1$ such
that $\bP(\log H_{m_1} \ge a(m_1)) > \e/2$. Inductively pick
$m_2,m_3,\ldots$ as follows. Given $m_1,m_2,\ldots,m_{k-1}$, let
$n_{k-1} = m_1+\cdots+m_{k-1}$. For sufficiently large $m_k\ge1$, we have
\[
\frac{\log H_{n_{k-1}}}{a(n_k)} \ge-\delta/2,\qquad\frac
{a(m_k)}{a(n_k)} \ge1-\delta/2\quad\mbox{and}\quad\bP\bigl(\log H'_{m_k} \ge
a(m_k) \mid \kS_0^{n_{k-1}}\bigr) > \e/2.
\]
Here, $n_k = n_{k-1} + m_k$ and $H'_{m_k} = H_{m_k}\circ
T_{(n_{k-1}/d,\ldots,n_{k-1}/d)}$. Note that such an $m_k$ always
exists, but depends on $n_{k-1}$ and $\log H_{n_{k-1}}$, so it is a
$\kS_0^{n_{k-1}}$-measurable random integer. Now, observe that
\[
\frac{\log H_{n_k}}{a(n_k)} \ge\frac{\log H_{n_{k-1}}}{a(n_k)} + \frac{\log H'_{m_k}}{a(n_k)} \ge- \delta/2 +
(1-\delta/2)\frac
{\log H'_{m_k}}{a(m_k)} \ge1 - \delta
\]
if $\log H'_{m_k} \ge a(m_k)$. Therefore,
\[
\bP \biggl(\frac{\log H_{n_k}}{a(n_k)} \ge1 - \delta \Big| \kS_0^{n_{k-1}}
\biggr) \ge\bP\bigl(\log H'_{m_k} \ge a(m_k)
\mid \kS _0^{n_{k-1}}\bigr) > \e/2.
\]
By L\'evy's extension of the second Borel--Cantelli lemma (see \cite{Wil1991}, page~124),
\[
\bP \biggl(\frac{\log H_{n_k}}{a(n_k)}\ge1 - \delta \mbox { i.o.} \biggr) = 1.
\]
Since $\delta>0$ is arbitrary, this gives (\ref{borneo}). In
particular, for $\lambda= \Lambda_q(V)$,
\[
\bP\bigl(\bar h_\infty^\lambda= \infty\bigr) = \bP \Bigl(\limsup
_{n\to
\infty}\log h_n^\lambda= \infty \Bigr) \ge\bP
\Bigl(\limsup_{n\to\infty}\log H_n = \infty \Bigr) = 1.
\]\upqed
\end{pf*}

\section*{Acknowledgments}
We thank F. Comets for providing us with an overview of open problems
regarding disorder regimes of directed polymers. A. Yilmaz thanks I. Corwin
and F. Rezakhanlou for valuable discussions. 

F.~Rassoul-Agha was partially supported by National Science Foundation
Grant DMS-14-07574 and the Simons Foundation grant 306576.
T.~Sepp\"al\"ainen was partially supported by National Science
Foundation Grant DMS-13-06777 and by the Wisconsin Alumni Research Foundation.
A.~Yilmaz was partially supported by European Union FP7 Grant
PCIG11-GA-2012-322078.



\printhistory
\end{document}